\documentclass{amsart}
\usepackage{hyperref}

\newtheorem{thm}{Theorem}[section]
\newtheorem{lemma}[thm]{Lemma}
\newtheorem{coro}[thm]{Corollary}
\newtheorem{hypo}[thm]{Hypothesis {\bf H.}\hspace*{-0.6ex}}
\newtheorem{rem}[thm]{Remark}

\newcommand{\bea}{\begin{eqnarray}}
\newcommand{\eea}{\end{eqnarray}}
\newcommand{\ba}{\begin{array}}
\newcommand{\ea}{\end{array}}

\newcommand{\bs}{\backslash}

\newcommand{\sig}{\sigma}
\newcommand{\lam}{\lambda}
\newcommand{\gam}{\gamma}

\newcommand{\bth}{\begin{thm}}
\newcommand{\eth}{\end{thm}}
\newcommand{\bl}{\begin{lemma}}
\newcommand{\el}{\end{lemma}}
\newcommand{\bk}{\begin{coro}}
\newcommand{\ek}{\end{coro}}
\newcommand{\bh}{\begin{hypo}}
\newcommand{\eh}{\end{hypo}}
\newcommand{\br}{\begin{rem}}
\newcommand{\er}{\end{rem}}
\newcommand{\bpf}{\begin{proof}}
\newcommand{\epf}{\end{proof}}

\numberwithin{equation}{section}

\begin{document}

\title{Reconstructing Jacobi Matrices from Three Spectra}

\author{Johanna Michor}
\address{Institut f\"ur Mathematik\\
Strudlhofgasse 4\\ 1090 Wien\\ Austria\\ and International Erwin Schr\"odinger
Institute for Mathematical Physics, Boltzmanngasse 9\\ 1090 Wien\\ Austria}
\email{Johanna.Michor@esi.ac.at}

\author{Gerald Teschl}
\address{Institut f\"ur Mathematik\\
Strudlhofgasse 4\\ 1090 Wien\\ Austria\\ and International Erwin Schr\"odinger
Institute for Mathematical Physics, Boltzmanngasse 9\\ 1090 Wien\\ Austria}
\email{Gerald.Teschl@univie.ac.at}
\urladdr{http://www.mat.univie.ac.at/\~{}gerald/}



\maketitle

\begin{abstract}
Cut a Jacobi matrix into two pieces by removing the $n$-th column and $n$-th
row. We give necessary and sufficient conditions for the spectra of the
original matrix plus the spectra of the two submatrices to uniquely determine
the original matrix. Our result contains Hochstadt's theorem as a special case.
\end{abstract}

\section{Introduction}

The topic of this paper is inverse spectral theory for Jacobi matrices, that
is, matrices of the form
\begin{equation}
H = \left( \begin{array}{ccccc} 
b_1 & a_1 &  &  &     \\
a_1 & b_2 & a_2 & & \\
 & \ddots & \ddots & \ddots &  \\
 &  & a_{N - 2} & b_{N - 1} & a_{N - 1} \\
 &  &  & a_{N - 1} & b_N \\
\end{array} \right)
\end{equation}
This is an old problem closely related to the moment problem (see
\cite{simp} and the references therein), which has attracted considerable interest recently
(see, e.g., \cite{gsfj} and the references therein, \cite{gla}, \cite{gibs}, \cite{shieh}).
In this note we want to investigate the following question: Remove the
$n$-th row and the $n$-th column from $H$ and denote the resulting submatrices
by $H_-$ (from $b_1$ to $b_{n-1}$) respectively $H_+$ (from $b_{n+1}$ to $b_N$).
When do the spectra of these three matrices determine the original matrix $H$?
We will show that this is the case if and only if $H_-$ and $H_+$ have no
eigenvalues in common.

{}From a physical point of view such a model describes a chain of $N$
particles coupled via springs and fixed at both end points 
(see \cite{tjac}, Section~1.5). Determining the eigenfrequencies
of this system and the one obtained by keeping one particle fixed, one 
can uniquely reconstruct the masses and spring constants. Moreover,
these results can be applied to completely integrable systems, in particular the
Toda lattice (see e.g., \cite{tjac}).

\section{Main result}

To set the stage let us introduce some further notation. We denote the spectra
of the matrices introduced in the previous section by
\begin{equation}
\sig(H) = \{ \lam_j \}_{j = 1}^{N}, \quad
\sig(H_-) = \{ \mu^-_k \}_{k = 1}^{n - 1}, \quad
\sig(H_+) = \{ \mu^+_l \}_{l = 1}^{N - n}.
\end{equation}
Moreover, we denote by $( \mu_j )_{j=1}^{N-1}$ the ordered eigenvalues of
$H_-$ and $H_+$ (listing common eigenvalues twice) and recall
the well-known formula
\begin{equation} \label{gzn}
g(z,n) = \frac{\prod_{j=1}^{N-1} (z - \mu_j)}{\prod_{j=1}^{N-1} (z - \lam_j)} =
\frac{-1}{z-b_n + a_n^2 m_+(z,n) + a_{n-1}^2 m_-(z,n)},
\end{equation}
where $g(z,n)$ are the diagonal entries of the resolvent $(H-z)^{-1}$ and $m_\pm(z,n)$
are the Weyl $m$-functions corresponding to $H_-$ and $H_+$. The Weyl functions
$m_\pm(z,n)$ are Herglotz and hence have a representation of the following form
\begin{eqnarray}
m_-(z, n) & = & \sum_{k = 1}^{n - 1} \frac {\alpha_k^-}{\mu_k^- - z}, 
\qquad \alpha_k^- > 0, \quad \sum_{k = 1}^{n - 1}  \alpha_k^- =1,\\
m_+(z, n) & = & \sum_{l = 1}^{N - n} \frac {\alpha_l^+}{\mu_l^+ - z}, 
\qquad \alpha_l^+ > 0,  \quad \sum_{l = 1}^{N - n} \alpha_l^+=1.
\end{eqnarray}

With this notation our main result reads as follows

\begin{thm}
To each Jacobi matrix $H$ we can associate spectral data
\begin{equation}
\{ \lam_j \}_{j = 1}^{N}, \quad ( \mu_j, \sig_j )_{j = 1}^{N - 1},
\end{equation}
where $\sig_j = +1$ if $\mu_j \in \sig(H_+) \bs \sig(H_-)$, 
$\sig_j = -1$ if $\mu_j \in \sig(H_-) \bs \sig(H_+)$, and
\begin{equation}
\sig_j = \frac {a_n^2 \alpha_l^+ - a_{n - 1}^2 \alpha_k^-}
{a_n^2 \alpha_l^+ + a_{n - 1}^2 \alpha_k^-}
\end{equation}
if $\mu_j=\mu_k^- = \mu_l^+$.

Then these spectral data satisfy
\begin{enumerate}
\item[(i)] $\lambda_1 < \mu_1 \leq \lambda_2 \leq \mu_2 \leq \dots  < \lambda_N$,
\item[(ii)] $\sig_j=\sig_{j+1} \in(-1,1)$ if $\mu_j = \mu_{j+1}$ and
$\sig_j \in \{ \pm 1\}$ if $\mu_j \ne \mu_i$ for $i\ne j$
\end{enumerate}
and uniquely determine $H$. Conversely, for every given set of spectral data
satisfying $(i)$ and $(ii)$, there is a corresponding Jacobi matrix $H$.
\end{thm}

\begin{proof}
We first consider the case where $H_-$ and $H_+$ have no eigenvalues in common.
The interlacing property (i) is equivalent to the Herglotz property of $g(z,n)$.
Furthermore, the residues $\alpha^-_i$ can be computed from (\ref{gzn})
\begin{eqnarray}   \nonumber
\frac {\prod_{j = 1}^N (z - \lambda_j)}
{\prod_{k = 1}^{n - 1} (z - \mu_k^-) \prod_{l = 1}^{N - n} (z - \mu_l^+)}
& = & z - b_n - a_n^2 \sum_{l = 1}^{N - n} \frac {\alpha_l^+}{z - 
\mu_l^+} \\    \label{inserting our ansatz}
&  & -\, a_{n - 1}^2 \sum_{k = 1}^{n - 1} 
\frac {\alpha_k^-}{z - \mu_k^-}.
\end{eqnarray}
and are given
by $ \alpha_i^-=  a_{n - 1}^{-2} \beta_i^-$, where 
\begin{equation}  \label{reconstruct betaim}
\beta_i^- = -\, \frac
{\prod_{j = 1}^N (\mu_i^- - \lambda_j)}
{\prod_{l \neq i} (\mu_i^- - \mu_l^-) \prod_{l = 1}^{N - n} 
(\mu_i^- - \mu_l^+)}, \qquad
a_{n - 1}^2 = \sum_{i = 1}^{n - 1} \beta_i^-.
\end{equation}
Similarly, $\alpha_l^+ = a_n^{-2} \beta_l^+$, where
\begin{equation} \label{reconstruct betaip}
\beta_l^+ = -\, \frac {\prod_{j = 1}^N (\mu_l^+ - \lambda_j)}
{\prod_{k = 1}^{n - 1} (\mu_l^+ - \mu_k^-) \prod_{p \neq l} 
(\mu_l^+ - \mu_p^+)}, \qquad
a_n^2 = \sum_{l = 1}^{N - n} \beta_l^+. 
\end{equation}
Hence $m_\pm(z,n)$ are uniquely determined and thus $H_\pm$ by standard
results from the moment problem. The only remaining coefficient $b_n$
follows from the well-known trace formula
\begin{equation} \label{reconstruct bn}
b_n  =  \mathrm{tr}(H) -  \mathrm{tr}(H_-) -  \mathrm{tr}(H_+) = \sum_{j = 1}^N \lambda_j 
- \sum_{k = 1}^{n - 1} \mu_k^- - \sum_{l = 1}^{N - n} \mu_l^+.
\end{equation}

Conversely, suppose we have the spectral data given.
Then we can define $a_n$, $a_{n-1}$, $b_n$, $\alpha_k^-$, $\alpha_l^+$
as above. By (i), $\alpha_k^-$ and $\alpha_l^+$ are positive and hence give rise
to $H_\pm$. Together with $a_n$, $a_{n-1}$, $b_n$ we have thus defined a
Jacobi matrix $H$. By construction, the eigenvalues  $\mu_k^-$, $\mu_l^+$
are the right ones and also (\ref{gzn}) holds for $H$. Thus $\lam_j$ are the
eigenvalues of $H$, since they are the poles of $g(z,n)$.

Next we come to the general case where  $\mu_{j_0} = \mu_{k_0}^- = \mu_{l_0}^+\,\, ( = \lambda_{j_0})$ at least for one $j_0$. Now some factors in the
left hand side of (\ref{inserting our ansatz}) will cancel and we can no longer
compute $\beta_{k_0}^-$, $\beta_{l_0}^+$, but only
$\gam_{j_0}= \beta_{k_0}^- + \beta_{l_0}^+$. However, by definition of
$\sig_{j_0}$ we have
\begin{equation}  \label{distribute sum}
\beta_{k_0}^- = \frac {1 - \sigma_{j_0}}{2}\, \gamma_{j_0}, \quad \quad 
\beta_{l_0}^+ = \frac {1 + \sigma_{j_0}}{2}\, \gamma_{j_0}.
\end{equation}
Now we can proceed as before to see that $H$ is uniquely determined by the spectral data.

Conversely, we can also construct a matrix $H$ from given spectral data, but it
is no longer clear that $\lam_j$ is an eigenvalue of $H$ unless it is a pole
of $g(z,n)$. However, in the case $\lam_{j_0}= \mu_{k_0}^- = \mu_{l_0}^+$ we
can glue the eigenvectors of $H_-$ and $H_+$ to give an eigenvector
corresponding to $\lam_{j_0}$ of $H$.
\end{proof}

The special case where we remove the first row and the first column (in which case $H_-$ is not present) corresponds to Hochstadt's theorem \cite{hspec}. Similar
results for (quasi-)periodic Jacobi operators can be found in \cite{ttr}.


\begin{thebibliography}{XXXX}
\bibitem{gsfj} F. Gesztesy and B. Simon, {\em $m$-functions and inverse spectral
analysis for finite and semi-infinite Jacobi matrices}, J. Anal. Math. {\bf 73},
267--297 (1997).
\bibitem{gla} G. M. L. Gladwell, {\em On isospectral spring-mass systems}, Inverse
 Probl. {\bf 11}, 591--602 (1995).
\bibitem{gibs} P. C. Gibson, {\em Inverse spectral theory of finite Jacobi 
matrices}, to appear in Memoires of the Amer.\ Math.\ Soc..
\bibitem{hspec} H.~Hochstadt, {\em On the construction of a Jacobi matrix from
 spectral data}, Lin.~Algebra Appl. {\bf 8}, 435--446, 1974.
\bibitem{shieh} C.-T. Shieh, {\em Reconstruction of a Jacobi matrix with mixed data},
preprint.
\bibitem{simp} B. Simon, {\em The classical moment problem as a self-adjoint finite
difference operator}, Advances in Math. {\bf 137}, 82--203 (1998).
\bibitem{ttr} G. Teschl, {\em Trace Formulas and Inverse Spectral Theory for
Jacobi Operators}, Comm. Math. Phys. {\bf 196}, 175--202 (1998).
\bibitem{tjac} G. Teschl, {\em Jacobi Operators and Completely Integrable Nonlinear Lattices}, Math. Surv. and Mon. {\bf 72}, Amer. Math. Soc., Rhode Island, 2000.
\end{thebibliography}
\end{document}